\input AHTOH-E.STY

\let\documentclass\relax

\documentclass

\UDC{
512.74       
+ 512.543.72 
}


\title{%
The dimension of solution sets to
systems of equations in algebraic groups
}
\author{%
Anton A. Klyachko
\quad
Maria A. Ryabtseva
}
\address{
\myAddressW
\quad ryabtseva.mariya@gmail.com
}

\grantsFirst{\RFBR19-01-00591}

\abstract{%
The Gordon--Rodriguez-Villegas theorem says that, in a finite group, the
number of solutions to a system of coefficient-free equations is divisible
by the order of the group if the rank of the matrix composed of the
exponent sums of \hbox{$j$-th} unknown in $i$-th equation is less than the
number unknowns. We obtain analogues of this and similar facts for
algebraic groups. In particular, our results imply that the dimension of
each irreducible component of the variety of homomorphisms from a finitely
generated group with infinite abelianisation into an algebraic group $G$
is at least $\dim G$.
}

\s 1.
Introduction

\proclaim{Solomon's theorem} \rm [Solo69].
In any group, the
number of solutions to a system of coefficient-free equations is
divisible by the order of the group if the system has less equations than
unknowns.

Here as usual, an \emph{equation over a group $G$}
is an
expression of the form $v(x_1,\dots,x_m)=1$, where~$v$ is a word whose
letters are unknowns, their inverses, and
elements of~$G$ called \emph{coefficients}
(there are no coefficients though by the condition of Solomon's theorem).
In other words, the left-hand side of an equation is an element of the
free product $G*F(x_1,\dots,x_m)$ of $G$ and the free group
$F(x_1,\dots,x_m)$ of rank $m$ (where $m$ is the number of unknowns).

Solomon's theorem was generalised in different directions, see
[Isaa70],
[Stru95],
[AmV11],
[GRV12],
[KM14],
[KM17],
[BKV18],
and literature cited therein.
For example, the following fact was proved in [KM14].

\Theorem KM \rm[KM14].
The number of solutions to a system of equations
over a group
is divisible by the order of the
centraliser of the set of coefficients if the rank of the matrix
composed of the exponent sums of
the $j$-th unknown in the $i$-th equation is less than the
number unknowns.

For coefficient-free equations, this theorem was obtained earlier by
Gordon and Rodriguez-Villegas [GRV12]
(and, in this case, the number of solutions is divisible by the order
of the whole group).
For example, the number of solutions to the system of equations
$\left\{x^{100}y^{100}[x,y]^{777}=1,\;(xy)^{\the\year}=1\right\}$
is always divisible by the order of the group because this system
(although not covered by Solomon's theorem)
satisfies the
conditions of the Gordon--Rodriguez-Villegas theorem: the matrix
composed of
the exponent sums (called the \emph{matrix of the system of equations})
has the form%
~$\left(
{\scriptscriptstyle100\;100 \atop
\scriptscriptstyle\the\year\;\the\year}
\right)
$,
and its rank is one (while the number of unknowns is two).

Note that, in these divisibility theorems,
the group
is not assumed to be necessarily finite; the divisibility is always
understood in the sense of cardinal arithmetic: each infinite cardinal is
divisible by all smaller nonzero cardinals. However, the most interesting
applications of these theorems concern finite
groups.

The purpose of this paper is to obtain analogues of the divisibility
theorems for equations over algebraic groups. An analogue of the number of
solutions to a system of equations is, in this case, the dimension
of the variety of solutions. Slightly more precisely one can say that
the dimension is something like the logarithm of the number of solutions.
Note that the solutions of a (finite or infinite) system of equations in
$m$ unknowns over an affine algebraic group $G$ form an affine
algebraic subvariety in $G^m$. An analogue of the Solomon theorem is the
following simple observation.

\Th 0.
In any affine algebraic group $G$, the
dimension of each irreducible component of the variety of solutions to a
finite
system of
equations \(possibly with coefficients\) is at least
$$
\Bigl((the\ number\ of\ unknowns)-
(the\ number\ of\ equations)\Bigr)\!\cdot\dim G.
$$

\Proof
The left-hand sides of $n$ equations
in $m$ unknowns
specify a morphism of algebraic
varieties $\gamma\:G^m\to G^n$. The solution variety is the
fiber~$\gamma^{-1}\bigl(1,\dots,1\bigr)$. Each irreducible component
$M$ of the variety $G^m$ is isomorphic to $G_0^m$ as an algebraic variety
(where~$G_0$~is the identity component of~$G$) and,
therefore, $\dim M=m\cdot\dim G$. The restriction of~$\gamma$
to~$M$ is a dominant morphism of irreducible algebraic varieties $M\to N$,
where~$N$ is the Zariski closure of~$\gamma(M)$ (this morphism is dominant 
because any set is dense in its closure).
It remains to refer to the following general fact (see, e.g., [Bo72]).

\proclaim Fiber lemma.
If $\gamma\:M\to N$ is a
dominant
morphism of irreducible algebraic varieties,
then, for any point~$s\in f(M)$,
the dimension of the fiber~$\gamma^{-1}(s)$ is at least
$\dim M-\dim N$; moreover, $N$ contains a nonempty open subset
$U \subset f(X)$ such that $\dim \gamma^{-1}(u)=\dim M-\dim N$
for any $u\in U$.

Theorem 0 is hardly new; for coefficient-free equations,
this fact was noted, e.g., in [LM11].

Solomon's theorem looks better than Theorem 0 in the sense that
Solomon's theorem asserts divisibility, while Theorem 0 only asserts an
inequality (this is unavailable --- the
divisibility of dimensions would be too much of course). In all other
aspects, Theorem 0 looks better; but one cannot improve Solomon's theorem
accordingly as the following simple examples show.

\proclaim False theorem.
In any group $G$, the
number of solutions to a consistent system of equations
with more unknowns than equations
\item{\rm1)}
is at least the order of the group;
\item{\rm2)}
is a multiple of
$|G|^{(the\ number\ of\ unknowns)-(the\ number\ of\ equations)}$
if the system is coefficient-free.

\Proof
\item{\rm1)}
The first assertion is false for a simple reason: in the symmetric group
of order six, the equation $x^3y^3=(1\;2\;3)$ has precisely three
solutions, because $x$ and $y$ must have the same parity. Moreover, they
cannot be even, because the
cube of an even permutation is identity. The product of two
different transpositions is always a 3-cycle; three of these
six products give one 3-cycle, and three give the other. (Surely,
the number of solutions must be a multiple of three by Theorem KM.)

\item{\rm2)}
{
The second assertion is also false as the equation $x^2y^2z^2=1$
in the same symmetric group shows.
The squares of permutations are even; therefore, for these squares,
we have three possibilities:
\itemitem{-}
they are the same 3-cycle;
\itemitem{-}
or
one of the squares is the identity permutation,
and the other two are different 3-cycles;
\itemitem{-}
or
all three square are the identity permutation.

\enditem

The
identity permutation has four square roots, a 3-cycle has one square
root, and we obtain $2+3\cdot2\cdot4+4\cdot4\cdot4=90$
solutions.
}

\medskip
\enditem
The following fact is an analogue of Theorem KM.

\Th 1.
If the rank of the matrix of a
system of equations
with finitely many unknowns
over an
affine algebraic group
is less than the number of unknowns,
then the
dimension of each component of the variety of solutions
to this system is at least
the dimension of the centraliser of the set of coefficients.

The equation $[x,y]=1$ (with zero matrix and two unknowns)
over a nonabelian connected group
shows that the estimate from Theorem 1 (even for coefficient-free case)
cannot be strengthened to the inequality
$$
\dim(the\ variety\ of\ solutions)\ge
\Bigl((the\ number\ of\ unknowns)-
(rank\ of\ the\ matrix\ of\ the\ system)\Bigr)\!\cdot\dim G
$$
(similar to Theorem 0).


\Question \rm(suggested by the referee).
Is it true that, under the conditions of Theorem 1,
the class of the solution variety is a multiple of the 
class of the centraliser of the coefficients in the Grothendieck
ring of varieties?  

If yes, then it would include both Theorem 1 and Theorem KM (for finite 
groups).  


\Corollary \rm(or rather a reformulation of Theorem 1 for
coefficient-free equations).
Suppose that $G$ is an affine algebraic
group and~$F$ is a finitely generated group whose commutator subgroup
is of infinite index. Then the dimension of each irreducible component
of the variety of homomorphisms $F\to G$ is at least $\dim G$.

The \emph{representation variety} (or
\emph{homomorphism variety})
$\Hom(F,G)$ of a finitely generated group
$F$ in an algebraic group $G$ was studied in many papers, see,
e.g.,
[RBCh96],
[MO10],
[LM11],
[LL13],
[LS05],
[Ki18],
[LT18]
and literature cited therein.
In particular, it was studied when there exist
injective or
\emph{topologically surjective} homomorphisms, i.e. homomorphisms with
a dense (in the Zariski topology) image, and how many such
homomorphisms are there. For example, in [Ki18] it was shown that
topologically surjective homomorphisms from the fundamental group of a
closed oriented surface of genus larger than one into a real semisimple
algebraic group always exist, and there are a lot of them in some sense.
In [BGGT12], it was shown that there exists a homomorphism from any free
group into almost any semisimple algebraic group such that the
restriction of the homomorphism to any noncyclic subgroup is topologically
surjective. See also a recent survey [GKP18].

Generally, neither injective nor
topologically surjective homomorphisms
form
a variety; however,
the following theorem shows that,
in some sense, all
``components" of these (and similar) non-varieties have dimension not
less than that of the commutator subgroup of $G$. We
consider the following property of a homomorphism $\phi\:F\to G$ from
a finitely generated group $F$ containing a subgroup $W$ into an
algebraic group $G$ containing a closed subgroup $A$.

{\parindent 2cm
\item{${\bf F}_{W,A}$} (faithfulness):
$\phi(W)\subseteq A$ and
the restriction of $\phi$ to $W\subseteq F$ is injective;
\item{${\bf S}_{W,A}$} (topological surjectivity):
the closure of the subgroup $\phi(W)$ is $A$;
\item{${\bf C}_{W,A}$} (connectivity):
$\phi(W)\subseteq A$ and
the closure of $\phi(W)$ is irreducible.

}

\Th 2.
Suppose that $A$ is a subgroup of an affine algebraic group $G$,
and $F$~is a finitely generated group whose commutator subgroup
is of infinite index.
Then each homomorphism $\alpha\:F\to G$ having some
(possibly infinite)
combination (conjunction)
$$
{\bf P}=
\(\bigwedge_{W\in\cal F}{\bf F}_{W,A}\)
\wedge
\(\bigwedge_{W\in\cal S}{\bf S}_{W,A}\)
\wedge
\(\bigwedge_{W\in\cal C}{\bf C}_{W,A}\),
\qbox{where $\cal F,S,C$ are some family of subgroups of $F$,}
$$
of the properties listed above 
is contained in a $(\dim[A,A])$-dimensional
irreducible subvariety of $\Hom(F,G)$ consisting of
homomorphisms having the same combination of properties $\bf P$
provided the index of the subgroup
$[F,F]\cdot\!\!\!\!\displaystyle\prod\limits_{W\in\cal F\cup C}\!\!\!\!W$
in~$F$ is infinite \rm(and no
conditions on the family $\cal S$ are imposed).

Here, one cannot replace
the dimension of the commutator subgroup with the dimension of
the entire group
as the following simple example shows:
there are many
topologically
surjective homomorphisms
from the infinite cyclic
group to the torus $(\C^*)^{\the\year}$
but they are all lone in the sense that any
non-zero-dimensional subvariety of $\Hom\(\Z,(\C^*)^{\the\year}\)$
contains a homomorphism that is not
topologically surjective (because any infinite
subvariety in $\C^*$ contains a finite-order element).

One cannot also remove the condition on the family $\cal F$ and $\cal C$
as simple examples of homomorphisms from $\Z$
to~$\SL_n(\C)$ show.

For instance, Theorem 2 implies that
\disp{\sl
if a homomorphism from the fundamental group
$F=\pres<x,y,z,t|[x,y]=[z,t]>$
of the genus-two orientable surface
{\rm(or from any other group with more
generators than relators)} into an algebraic group $G$
is
injective on the subgroup $\gp{x,y,z}$
and maps any non-abelian subgroup onto a
subgroup dense in $G$ , then it is contained in a $(\dim[G,G])$-dimensional
subvariety of $\Hom(F,G)$ consisting of homomorphisms with
the same two properties.
}
Partially, Theorem 2 is an analogue of the following known fact
[KM17]:
\disp{\sl
the number of surjective homomorphisms from
a finitely generated group whose commutator subgroup
is of infinite index into a {\rm(finite if you like)}
group $G$ is divisible by the order of the commutator subgroup
of~$G$.
{\rm For example,} the number of pairs of elements generating
a
\(2-generated\)
group
is always divisible by the order of the commutator subgroup
of this group.
}%
However, in ``algebraic settings"
some new effects arise.

Theorems 1 and 2 are special cases of a more general \emph{Main theorem},
which is stated in the next section. Our main theorem
is a direct ``algebraic" analogue of the main theorem
from [KM17] (whose special cases are
all
theorems on divisibility mentioned above
and some other curious facts, see [KM17]
and [BKV18]).

The authors thank D. A. Timashev and an anonymous referee for valuable 
remarks.

{\noindent \bf Our Notation and conventions}
are mainly standard. Note only that, if
$k\in \Z$, and $x$ and $y$ are elements some group, then $x^y$, $x^{ky}$,
and $x^{-y}$ denote $y^{-1}xy$, $y^{-1}x^ky$, and $y^{-1}x^{-1}y$,
respectively.
The commutator subgroup of a group~$G$ is denoted by $G'$
or $[G,G]$.
If $X$~is a
subset of a group, then
$\gp X$, $\nc X$, and $C(X)$
mean
the subgroup generated by~$X$, normal closure of~$X$, and the centraliser
of~$X$,
respectively,
The index of a subgroup $H$ of a group
$G$ is denoted $|G:H|$. The letter~$\Z$ denotes
the set of integers.
The symbol~$\Hom(A,B)$ denotes the set of homomorphisms from a group~$A$ to
a group $B$.
The free group of rank~$n$ is denoted as
$F(x_1,\dots,x_n)$ or $F_n$.
The symbol~$A*B$ denotes the free product of groups $A$ and~$B$.
The word \emph{variety} always means a 
(not necessarily irreducible)
quasi-projective algebraic
variety over an algebraically closed field
(of arbitrary characteristic); the word \emph{subvariety} means a locally
closed subset of a variety. All topological terms refer to the
Zariski topology.

\s 2.
Main theorem

A group $F$ equipped with an epimorphism $F\to\Z$ is called an
\emph{indexed} group. This epimorphism $F\to\Z$ is called the
\emph{degree} or \emph{indexing}
and denoted as $\deg$; thus, for any element $f$ of an
indexed group $F$, there is an integer $\deg f$, the group
$F$ contains elements of all integer degrees, and
${\deg(fg)=\deg f+\deg g}$ for any $f,g\in F$.

Suppose that $\phi\:F\to G$ is a homomorphism from an indexed group $F$ to
some group $G$, and $H$ is a subgroup of~$G$. We call the subgroup
$$
H_\phi=\bigcap_{f\in F}H^{\phi(f)}\cap C(\{\phi(f)\;|\;\deg f=0\})
$$
the \emph{$\phi$-core} of $H$. In other words,
the $\phi$-core $H_\phi$ of~$H$ consists of its
elements~$h$ such that $h^{\phi(f)}\in H$ for all $f$ and
$h^{\phi(f)}=h$ if $\deg f=0$.

\proclaim{Main theorem}.
Suppose that $H$ is an affine
algebraic
subgroup of an
algebraic
group $G$, and a
subvariety $\Phi$ of~$\Hom(F,G)$,
where $F$ is a
finitely generated
indexed group, has
the following two properties.
\item{\rm I.}
$\Phi$ is invariant under conjugation by elements of $H$:
$$
\qbox{if $h\in H$ and $\phi\in\Phi$, then the homomorphism
$\psi\:f\mapsto\phi(f)^h$ also lies in $\Phi$.}
$$
\item{\rm II.}
For any $\phi\in\Phi$ and any element $h$ of the $\phi$-core
$H_\phi$ of~$H$, the homomorphism $\psi$ defined as
$$
\psi(f)=
\cases{
\phi(f)& for all elements $f\in F$ of degree zero;
\cr
\phi(f)h& for
some element
$f\in F$ of degree one
\small(and, hence, for all elements of degree one)
\cr
}
$$
\nobreak
belongs to $\Phi$ too.
\enditem
Then the dimension of each irreducible component of~$\Phi$
is at least
the dimension of~$H$.

\goodbreak

Note that the mapping $\psi$ from Condition I is a homomorphism
for any $h\in G$, and the formula for $\psi$ from condition II defines
a
homomorphism for any $h\in C(\phi(\ker\deg))$ ([KM17], Lemma 0).
Conditions I and II only require these homomorphisms to belong to $\Phi$
(under some additional restrictions on $h$).

This theorem is an analogue of the main theorem in [KM17],
which asserts that, for any (abstract) group $G$ and any its
subgroup $H$, the cardinality of a set $\Phi$ of homomorphisms from
an indexed group $F$ to $G$ is divisible by~$|H|$ if
Conditions I and II hold.

\s 3.
Proof of Theorem 1

Let $A\subseteq G$ be the subgroup generated by all coefficients of all
equations, and let $F$ be the quotient group
$$
F=(A*F(x_1,\dots,x_n))/\nc{\{v_i\}}
$$
of the free product
$A*F(x_1,\dots,x_n)$ of $A$ and the free group
$F(x_1,\dots,x_n)$
by the normal subgroup
$\nc{\{v_i\}}$ generated by the left-hand sides of all equations. As
the set $\Phi$ we consider
the set of
homomorphisms $F\to G$ which restrict to the identity on $A$ (we assume
that the natural mapping $A\to F$ is an embedding
because otherwise equations have no solutions and we have
nothing to prove). Clearly, the solutions to the system of equations are
in a natural one-to-one correspondence with the elements of $\Phi$, and
$\Phi$
is a subvariety of $\Hom(F,G)$.

The condition on the rank means that $F$ admits an epimorphism onto $\Z$
whose kernel contains $A$. Now, take
the centraliser of $A$ in $G$ as $H$. The conditions of the main theorem
are obviously fulfilled. Indeed, Condition I holds because $h$
centralises $A\subseteq G$ and, hence, $\psi$ coincides with $\phi$ on
$A\subset F$; Condition II holds because elements of $A\subset F$ have
degree zero and, therefore, $\psi$ coincides with $\phi$ on
$A\subset F$ --- again.

\s 4.
Proof of Theorem 2

If the commutator subgroup of a finitely generated group is
of infinite index,
then, as is known, the group admits an
epimorphism onto $\Z$.
Therefore, the conditions of Theorem~2 implies the existence of an
indexing $\deg\:F\to\Z$ whose kernel contains all subgroup from the
families $\cal C$ and $\cal F$ (and the commutator subgroup $F'$ of
course). Let $H$ be the commutator subgroup $A'$ of $A$, and let $\Phi$ be
the set of homomorphisms $F\to G$ coinciding with a given
homomorphism~$\alpha\:F\to G$ modulo~$H$ and coinciding with $\alpha$ on
elements of
degree zero:
$$
\Phi=
\{\phi\:F\to G\;|\;\quad \phi(f)H=\alpha(f)H \hbox{ for all $f\in F$};
\quad
\phi(f)=\alpha(f) \hbox{ for all $f\in\ker\deg$} \}.
$$
Clearly, the conditions of the main theorem hold for these
$F$, $\deg$, $\Phi$, and $H$. Therefore,
the dimension of each component of the variety
$\Phi$ is at least
$\dim H=\dim A'$.
It remains to verify that, if the homomorphism $\alpha$
has the property $\bf P$,
then all homomorphisms from $\Phi$
have this property too. Suppose that $\phi\in\Phi$, i.e.
$\phi(f)=\alpha(f)h_f$ for all~$f\in F$, where $h_f\in H=A'$.

First, note that the subgroup $\phi(W')=\alpha(W')$ is dense in~$H=A'$
for all $W\in\cal S$.  Indeed, any algebraic group has finite
\emph{commutator width}, i.e. the morphism $\kappa\:A^{2n}\to A'=H$
sending
$(x_1,\dots,x_n,y_1,\dots,y_n)$ to the product of
commutators~$\prod[x_i,y_i]$ is surjective for sufficiently large $n$
(see, e.g., [VO88]). The image of a dense set under a continuous
surjective mapping is dense (and $D^{2n}$ is dense in $A^{2n}$ if
$D$ is dense in $A$). Thus, any nonempty open subset of~$A'=H$
contains an element of
$\kappa\(\bigl(\alpha(W)\bigr)^{2n}\)\subseteq\alpha(W')=\phi(W')$
as required.

Let $U\subseteq A$ be a nonempty open set.
Since $\alpha(W)$ is dense in $A$ (where $W\in\cal S$),
the
set~$U$ contains
an element~$\alpha(w)=\phi(w)h_w^{-1}$,
where $w\in W$.
Therefore, the open set
$\phi(w^{-1})U$ contains
$h_w^{-1}\in H=A'$.
Hence,
$\phi(w^{-1})U\cap A'$ is a nonempty open subset in~$A'=H$. As shown
above, it contains some element $\phi(w_1)$,
where $w_1\in W$.
Therefore, $U\ni\phi(ww_1)$,
and this completes the proof of property~${\bf S}_{W,A}$ for~$\phi$.

Properties ${\bf F}_{W,A}$ (where $W\in\cal F$) and ${\bf C}_{W,A}$
(where $W\in\cal C$) hold for the obvious reason:
all these subgroup $W$ consist of elements of degree zero
(by the choice
of~$\deg$), therefore, $\phi$ and $\alpha$ coincide
on all such subgroups $W$.

\s 5.
Proof of the main theorem

We follow the proof
of the main theorem
in~[KM17] with necessary modifications.
The first difficulty is that the subgroup
$\ker\deg\subset F$ may
be non-finitely generated and, therefore,
the set of homomorphisms
$\ker\deg\to G$ may have no
natural structure of an
algebraic variety.
To overcome this unpleasant feature we use
the following simple observation:
\disp{\sl
there exists a finite subset
$K\subseteq\ker\deg\subset F$ such that any two homomorphism
from $F$ to $G$ coinciding on $K$ coincide on $\ker\deg$.
}%
Indeed, suppose that $\ker\deg=\{d_1,d_2,\dots\}$ and
let
$\Pi_i$ be the set of pairs of homomorphisms $F\to G$ coinciding
on~$d_1,\dots,d_i$. Clearly, $\Pi_i$ form a decreasing chain of
subvarieties in $\Hom(F,G)\times\Hom(F,G)$. Such a chain must
stabilise: $\Pi_n=\Pi_{n+1}=\dots$ for some $n$.
Therefore, we can put $K=\{d_1,d_2,\dots,d_n\}$.

\medskip

We also need a similar fact:
\disp{\sl
there exists a finite subset
$A\subset F$ such that, if, for two homomorphisms $\alpha$
and $\beta$
from~$F$ to $G$, their
composition with the natural mapping~$G\to G/H$
\(where $G/H$ is the set of left cosets of $H$ in $G$ \)
coincide on $A$, then they
coincide on the entire group $F$:
}%
\vskip-3mm
$$
\forall \alpha,\beta\in\Hom(F,G)
\qquad
\Bigl(\forall a\in A\quad\alpha(a)H=\beta(a)H\Bigr)
\imp
\Bigl(\forall f\in F\quad\alpha(f)H=\beta(f)H\Bigr).
$$
Indeed, suppose now that
$F=\{d_1,d_2,\dots\}$ and $\Pi_i$ is the set of pairs
of homomorphisms
$\alpha,\beta\:F\to G$ such that
$\alpha(d_k)H=\beta(d_k)H$ for $k\le i$. Clearly,
$\Pi_i$ form a decreasing chain of subvarieties. Such a chain must
stabilise: $\Pi_n=\Pi_{n+1}=\dots$ for some $n$.
Therefore, we can put $A=\{d_1,d_2,\dots,d_n\}$.

\medskip

Now, consider the variety $X=G^K\times(G/H)^A$
consisting of all pairs of mappings $K\to G$ and~$A\to G/H$
(where $G/H$ is the variety of left cosets the subgroup
$H$ in~$G$). The group $H$ acts on~$X$ (by conjugations):
$$
h\circ(\alpha,\beta)=
\bigl(f\mapsto h\alpha(f)h^{-1},\; a\mapsto h\beta(a)\bigr).
$$

The \emph{tail} $\chi(\phi)$ of a
homomorphism $\phi\:F\to G$ is the
pair $(\phi_0,\phi_H)$, where $\phi_0$ is restriction of~$\phi$
to the set $K\subseteq\ker\deg\subset F$, and
$\phi_H\:A\to\{gH\;;\;g\in G\}$ is the mapping from $A$ to
$G/H$ which maps $a\in A$ to the coset~$\phi(a)H$.
Clearly, the mapping $\chi\:\Hom(F,G)\to X$ is a morphism of
algebraic varieties.

We say that two homomorphism $\phi,\psi\in\Phi$ are \emph{similar} and
write $\phi\sim\psi$ if their tails lie in the same orbit under the
action of
$H$ on $X$ described above.
Note that neither the similarity of homomorphisms nor coincidence
of their tails depend on the choice of sets $A$ and $K$
(these sets are needed solely
to make ``taking the tail" a morphism of
algebraic varieties and to make
the action of $H$ on tails an action of an
algebraic group on an algebraic variety):
$$
\eqalignno{
\chi(\phi)=\chi(\psi)\ \iff
&\cases{
\psi(f)=\phi(f)
&
for all $f\in F$ of degree zero and
\cr
\psi(f)H=\phi(f)H
&
for all $f\in F$.
}
&(*)
\cr\cr
\phi\sim\psi\ \iff
\hbox{ for some $h\in H$ }
&\cases{
\psi(f)=h\phi(f)h^{-1}&
for all $f\in F$ of degree zero and
\cr
\psi(f)H=h\phi(f)H&
for all $f\in F$.
}
}
$$



Without loss generality, we assume that the group
$H$ is irreducible (because the identity component
of~$H$ is a group of the same dimension).

Each class of similar homomorphisms is a locally closed
subvariety in~$\Hom(F,G)$ because
this class
is the preimage of an orbit under the morphism $\chi$, and the orbit of
an action of an algebraic group on an algebraic variety is always locally
closed (see, e.g., [VO88]). The main theorem
follows immediately from the
following proposition.

\Proposition.
The dimension of
each component of
each class of similar homomorphisms in
$\Phi$ is $\dim H$.
More precisely, for each $\phi\in\Phi$,
\item{\rm 1)}
the dimension of the
variety $X_\phi$ of
tails of homomorphisms from $\Phi$ similar to $\phi$
equals $\dim H-\dim H_\phi$;
\item{\rm 2)}
for each homomorphism $\psi$ similar to $\phi$, the dimension
of
each component of the variety of
homomorphisms from~$\Phi$ with the same
tail as $\psi$
equals~$\dim H_\phi$.

\Proof
To prove Assertion 1) note that the set
$\chi(\Phi)\subseteq X$ is invariant with respect to the
action of $H$ on $X$. Indeed,
$h\circ\chi(\phi)=\chi(f\mapsto h\phi(f)h^{-1})$.
This homomorphism $f\mapsto h\phi(f)h^{-1}$
lies in $\Phi$ by Condition I of the main theorem.
The tails of homomorphisms similar to $\phi$ is the orbit of the tail
of $\phi$ under this action. The dimension of an orbit is, as is
known, equal to the codimension of the stabiliser (this is a special case
of the fiber lemma). It remains to note that the
subgroup~$H_\phi$ is the stabiliser of the tail of~$\phi$ (by
formula~$(*)$).

Let us prove the second assertion. Choose an element $x\in F$ of degree
one. A homomorphism $\alpha\:F\to G$ is uniquely determined by its tail
and the value $\alpha(x)$ (by formula $(*)$). Moreover, for two
homomorphisms $\alpha$ and $\beta$ with the same tail, the quotient
$h=(\alpha(x))^{-1}\beta(x)$ must stabilise this tail, i.e.
lie in $H_\alpha$. Indeed, for all $f\in F$ of degree zero we have
\def\={\buildrel*\over=}
$$
\alpha(f^x)^h= \alpha(f)^{\alpha(x)h}=
\alpha(f)^{\beta(x)}\=\beta(f)^{\beta(x)}=\beta(f^x)\=\alpha(f^x),
\qbox{i.e. $h$ centralises the subgroup $\alpha(\ker\deg)$}
$$
and, for any element $f\in F$, we have
$$
\alpha(x)\alpha(f)H=
\alpha(xf)H\=
\beta(xf)H=
\beta(x)\beta(f)H=
\alpha(x)h\beta(f)H\=
\alpha(x)h\alpha(f)H,
\qbox{i.e. $h\in \alpha(f)H\alpha(f)^{-1}$.}
$$
(Here, equalities $\=$ follow from formula $(*)$.)
Thus, $h=(\alpha(x))^{-1}\beta(x)\in H_\alpha$.

On the other hand, if $h$ is an element of $H_\alpha$,
then the mapping
$$
f\mapsto\cases{
               \alpha(f),  &if $\deg f=0$\cr
               \alpha(x)h, &if $f=x$\cr
               }
\eqno{(**)}
$$
can obviously be extended to a
homomorphism with the same tail as $\alpha$
(see the remark after the statement of the main theorem in Section~2).
This homomorphism lies in $\Phi$ by Condition II of the main theorem.

We showed that,
for any $\alpha\in\Phi$, the
mapping $H_\alpha\to\Hom(F,G)$ that maps an element $h\in H_\alpha$
to the homomorphism $(**)$ is an injective morphism of algebraic
varieties whose image is the set of elements of
$\Phi$ with the same tail as $\alpha$. This means that
$\dim\chi^{-1}(\chi(\alpha))=\dim H_\alpha$.

It remains to note that, for similar homomorphisms $\psi$ and $\phi$,
the subgroup
$H_\phi$ and $H_\psi$ are isomorphic and even
conjugate in $H$ because they are the stabilisers of points
$\chi(\phi)$ and $\chi(\psi)$ lying in the same orbit under the
action of $H$ on  $X$.
This completes the proof of Assertion 2).

The variety $X_\phi$ is irreducible (because we
assume that the group $H$ is irreducible) and each
component $\Pi_\phi\ni\phi$
of the variety of homomorphisms similar to $\phi$ maps onto
$X_\phi$ surjectively because
$\chi\bigl(f\mapsto h\phi(f)h^{-1}\bigr)=h\circ\chi(\phi)$
and the variety
$\{f\mapsto h\phi(f)h^{-1}\;|\;h\in H\}$ is irreducible (as it is an orbit
under the action of the connected group $H$ on~$\Hom(F,G)$). Therefore,
1), 2), and
the following general fact imply that
the dimension of each component of each class of similar homomorphisms
from $\Phi$ is $\dim H$.

\Lemma.
Let $\gamma$ be a
morphism from a variety~$P$ to an irreducible variety $N$ such that
all components of all fibers~$\gamma^{-1}(y)$ \(where $y\in N$\)
have the same dimension $d$
and the image of each component of~$P$ is dense in $N$.
Then the dimension of each component of~$P$ is $\dim N+d$.

\Proof
Let $M$ be an irreducible component of~$P$. Take an open
nonempty subset $U\subseteq \gamma(M)\subseteq N$
from the fiber lemma.
Irreducibility of $M$ implies that
$\gamma^{-1}(U)$ contains a point $x$ not belonging to other
components of~$P$.
Thus, each component $K$ of the fiber $\gamma^{-1}(\gamma(x))$
containing $x$ must lie in $M$ and, therefore, be a
component of the fiber of the restriction of the morphism $\gamma$ to $M$.
Clearly, the other components of this
fiber~$M\cap\gamma^{-1}(\gamma(x))$ have dimension at most $d=\dim K$.
By the fiber lemma, we obtain $\dim M=\dim K+\dim N$ as
required. This completes the proofs of the lemma and main
theorem.

\baselineskip11.2pt

\References






[Bo72]
A. Borel,
Linear algebraic groups.
Benjamin, New York, 1969.

[VO88]
E. B. Vinberg, A. L. Onishchik,
Seminar on Lie groups and algebraic groups.
Springer Verlag, 1990.

[GKP18]
N. L. Gordeev, B. E. Kunyavski\u\i, E. B. Plotkin,
Geometry of word equations in simple algebraic groups
over special fields,
Russian Mathematical Surveys, 73:5(443) (2018), 753--796.
\arXiv 1808.02303

[Stru95]   %
S. P. Strunkov,
On the theory of equations in finite groups,
Izvestiya: Mathematics, 59:6 (1995), 1273-1282.

[AmV11]     %
A. Amit, U. Vishne,
Characters and solutions to equations in finite groups,
J. Algebra Its Appl., 10:4 (2011), 675-686.

[BGGT12]
E. Breuillard, B.Green, R. Guralnick, and T. Tao,
Strongly dense free subgroups of semisimple algebraic groups,
Israel Journal of Mathematics, 192:1 (2012), 347-379.
\arXiv 1010.4259



[BKV18]
E. K. Brusyanskaya, A. A. Klyachko, A. V. Vasil'ev,
What do Frobenius's, Solomon's, and Iwasaki's theorems
on divisibility in groups have in common?,
Pacific Journal of Mathematics, 302:2 (2019), 437-452.
\arXiv 1806.08870


[GRV12]  %
C. Gordon, F. Rodriguez-Villegas,
On the divisibility of $\#\Hom(\Gamma, G)$ by $|G|$,
J. Algebra, 350:1 (2012),
300-307.
\arXiv 1105.6066




[Isaa70]
I. M. Isaacs,
Systems of equations and generalized characters in groups,
Canad. J. Math., 22 (1970),
1040-1046.


[Ki18]
K. Kishore,
Representation variety of surface groups,
Proc. Amer. Math. Soc. 146 (2018), 953-959.
\arXiv 1702.05981

[KM14]
A. A. Klyachko, A. A. Mkrtchyan,
How many tuples of group elements have a given property?
With an appendix by Dmitrii V. Trushin,
Intern. J. of Algebra and Comp. 24:4 (2014), 413-428.
arXiv 1205.2824

[KM17]
A. A. Klyachko, A. A. Mkrtchyan,
Strange divisibility in groups and rings,
Arch. Math. 108:5 (2017), 441-451.
\arXiv 1506.08967



[LL13]
M. Larsen, A. Lubotzky,
Representation varieties of Fuchsian groups,
From Fourier analysis and number theory to radon transforms and
geometry, 375-397, Dev. Math., 28, Springer, New York, 2013.
\arXiv 1203.3408

[LS05]
M. Liebeck, A. Shalev,
Fuchsian groups, finite simple groups and representation varieties,
Inventiones mathematicae 159:2 (2005), 317-367.

[LM11]
S. Liriano, S. Majewicz,
Algebro-geometric invariants of groups
(the dimension sequence of representation variety),
Int. J. Algebra Comput., 21:4 (2011), 595-614.

[LT18]
D. D. Long, M. B. Thistlethwaite,
The dimension of the Hitchin component for triangle groups,
Geometriae Dedicata
(to appear).


[MO10]
J. Mart\'\i n-Morales, A. M. Oller-Marc\'en,
On the number of irreducible components of the representation
variety of a family of one-relator groups,
Internat. J. Algebra Comput. 20:1 (2010), 77-87.
\arXiv 0805.4716




[RBCh96]
A. S. Rapinchuk, V. V. Benyash-Krivetz, V. I. Chernousov,
Representation varieties of the fundamental groups
of compact orientable surfaces,
Israel Journal of Mathematics, 93:1 (1996), 29-71.



[Solo69]
L. Solomon,
The solutions of equations in groups,
Arch. Math., 20:3 (1969), 241-247.




\end